\documentclass[12pt,leqno]{article}
\usepackage{amsfonts}
\usepackage{amsmath, amsthm, amsfonts, amssymb, color}
\usepackage{mathrsfs}
\usepackage{color}
\theoremstyle{definition}

\newcommand{\scr}[1]{\mathscr #1}
\definecolor{wco}{rgb}{0.5,0.2,0.3}

\numberwithin{equation}{section} \theoremstyle{remark}

\newcommand{\ua}{\uparrow}

\title{{\bf Spectral Gap  for   Measure-Valued Diffusion Processes  }\footnote{Supported in
 part by  NNSFC (11771326, 11831014, 11726627).} }
\author{
{\bf Panpan Ren$^{b,c)}$,    Feng-Yu Wang$^{a,b)}$  }\\
\footnotesize{$^{a)}$ Center for Applied Mathematics, Tianjin University, Tianjin 300072, China}\\
 \footnotesize{$^{b)}$ Department of Mathematics,
Swansea University, Singleton Park, SA2 8PP, United Kingdom}\\
\footnotesize{$^{c)}$ Mathematical Institute,Woodstock Road, OX2 6GG, University of Oxford}\\
\footnotesize{ 673788@swansea.ac.uk, Panpan.ren@maths.ox.ac.uk wangfy@tju.edu.cn, F.-Y.Wang@swansea.ac.uk}}
\begin{document}
\allowdisplaybreaks
\def\R{\mathbb R}  \def\ff{\frac} \def\ss{\sqrt} \def\B{\mathbf
B}
\def\N{\mathbb N} \def\kk{\kappa} \def\m{{\bf m}}
\def\ee{\varepsilon}\def\ddd{D^*}
\def\dd{\delta} \def\DD{\Delta} \def\vv{\varepsilon} \def\rr{\rho}
\def\<{\langle} \def\>{\rangle} \def\GG{\Gamma} \def\gg{\gamma}
  \def\nn{\nabla} \def\pp{\partial} \def\E{\mathbb E}
\def\d{\text{\rm{d}}} \def\bb{\beta} \def\aa{\alpha} \def\D{\scr D}
  \def\si{\sigma} \def\ess{\text{\rm{ess}}}\def\s{{\bf s}}
\def\beg{\begin} \def\beq{\begin{equation}}  \def\F{\scr F}
\def\Ric{\mathcal Ric} \def\Hess{\text{\rm{Hess}}}
\def\e{\text{\rm{e}}} \def\ua{\underline a} \def\OO{\Omega}  \def\oo{\omega}
 \def\tt{\tilde}\def\[{\lfloor} \def\]{\rfloor}
\def\cut{\text{\rm{cut}}} \def\P{\mathbb P} \def\ifn{I_n(f^{\bigotimes n})}
\def\C{\scr C}      \def\aaa{\mathbf{r}}     \def\r{r}
\def\gap{\text{\rm{gap}}} \def\prr{\pi_{{\bf m},\varrho}}  \def\r{\mathbf r}
\def\Z{\mathbb Z} \def\vrr{\varrho} \def\ll{\lambda}
\def\L{\scr L}\def\Tt{\tt} \def\TT{\tt}\def\II{\mathbb I}
\def\i{{\rm in}}\def\Sect{{\rm Sect}}  \def\H{\mathbb H}
\def\M{\scr M}\def\Q{\mathbb Q} \def\texto{\text{o}} \def\LL{\Lambda}
\def\Rank{{\rm Rank}} \def\B{\scr B} \def\i{{\rm i}} \def\HR{\hat{\R}^d}
\def\to{\rightarrow}\def\l{\ell}\def\iint{\int}
\def\EE{\scr E}\def\Cut{{\rm Cut}}\def\W{\mathbb W}
\def\A{\scr A} \def\Lip{{\rm Lip}}\def\S{\mathbb S}
\def\BB{\scr B}\def\Ent{{\rm Ent}} \def\i{{\rm i}}\def\itparallel{{\it\parallel}}
\def\g{{\mathbf g}}\def\Sect{{\mathcal Sec}}\def\T{\mathcal T}\def\BB{{\bf B}}
\def\f{\mathbf f} \def\g{\mathbf g}\def\BL{{\bf L}}\def\MM{\mathbb M} \def\BG{{\mathbb G}}
\def\Bd{{\nn^{ext}}} \def\BdP{\nn^{ext}_\phi} \def\Bdd{{\bf \dd}} \def\Bs{{\bf s}} \def\GA{\scr A}
\def\Bg{{\bf g}}  \def\Bdd{{\bf d}} \def\supp{{\rm supp}}\def\div{{\rm div}}
\def\ddiv{{\rm div}}\def\osc{{\bf osc}}\def\1{{\bf 1}}\def\BD{\mathbb D}
\maketitle

\begin{abstract}  The spectral gap is estimated for some measure-valued processes, which are induced by the intrinsic/extrinsic derivatives on the space of finite measures over a Riemannian manifold. These processes are  symmetric with respect to the  Dirichlet and Gamma distributions arising from population genetics. In addition to the evolution of allelic frequencies investigated in the literature, they also describe stochastic movements of individuals.   \end{abstract} \noindent
 \vskip 2cm

\section{Introduction}

The Dirichlet  distribution  arises naturally in Bayesian inference as conjugate priors for categorical distribution and infinite  non-parametric discrete distributions respectively.   In population genetics, it describes   the distribution of  allelic frequencies (see for instance \cite{ConMoi69,J1,M2}).   To simulate the Dirichlet distribution using   stochastic dynamic systems,  some diffusion processes generalized from  the Wright-Fisher
diffusion have been considered, see  \cite{QQ,FMW,FW07,FW14,STA} and  references within.    In this paper, we investigate  diffusion processes induced by  the Dirichlet distribution and
the intrinsic/extrinsic   derivatives, where the extrinsic derivative term determines  the evolution of  allelic frequencies, and the intrinsic derivative drives the movement of individuals.    

In the following three subsections, we introduce  the reference measures, intrinsic and extrinsic derivatives, and the main results of the paper respectively. We will take the notation $\mu(f)=\int_E f\d\mu$ for
a measurable space $(E,\scr B,\mu)$ and $f\in L^1(E,\mu)$.

\subsection{Reference measures}
Let $(M,\<\cdot,\cdot\>_M)$ be a complete  Riemannian manifold. Consider    the space $\MM$ of all nonnegative finite measures on $M$,
and let $\MM_1:=\{\mu\in \MM: \mu(M)=1\}$ be the set of all probability measures on $M$. According to \cite[Theorem 3.2]{MR}, both spaces are Polish under the weak topology.
In general, for an ergodic Markov process $X_t$ with stationary distribution $\mu$, one may simulate $\mu$ by using the empirical measures $\mu_t:=\ff 1 t \int_0^t \dd_{X_s}\d s,$ where $\dd_{X_s}$ is the Dirac measure at point $X_s$. In practice, one may also approximate $\mu$ using the discrete time empirical measures
$$\tt\mu_n:=\ff 1 n\sum_{i=1}^n \dd_{X_i},\ \ n\ge 1.$$
See for instance \cite{Gao} and references within for the study of the convergence rate.

For  $0\ne \theta\in \MM$, the Dirichlet distribution $\BD_\theta$ with shape $\theta$ is the unique probability measure on $\MM_1$ such that for any measurable partition $\{A_i\}_{1=1}^n$ of $M$,
$$\MM_1\ni\mu\mapsto (\mu(A_1),\cdots,\mu(A_n))$$ obeys the Dirichlet distribution with parameter  $(\theta(A_1),\cdots,\theta(A_n))$. Recall that for any $0\ne \aa=(\aa_1,\cdots,\aa_n)\in [0,\infty)^n$, the Dirichlet distribution with parameter $\aa$ is the following  probability measure on the simplex $\{s=(s_1,\cdots,s_n): s_i\ge 0, \sum_{i=1}^n s_i=1\}$:
$$\mathbb D_\aa(\d s_1,\cdots, \d s_n):= \ff{\GG(\aa_1+\cdots+\aa_n)}{\GG(\aa_1)\cdots\GG(\aa_n)}  s_n^{\aa_n-1} \dd_{1-\sum_{1\le i\le n-1}s_i}(\d s_n)  \prod_{i=1}^{n-1} s_i^{\aa_i-1} \d s_i,$$   where in case $\aa_i=0$ we set
$\ff 1 {\GG(0)} s_i^{-1} \d s_i=\dd_0$, and $\dd_x$ denotes   the Dirac measure at point $x$ in a measurable space.
If  $\BD_\theta$ refers to the distribution of population on $M$, then under the state $\mu\in \MM_1$,  $\mu(A_1),\cdots,\mu(A_n)$ stand for the  proportions of population located in the areas $A_1,\cdots, A_n$ respectively.

 We will also consider the Gamma distribution $\BG_\theta$ on $\MM$ with shape $\theta$, whose marginal distribution on $\MM_1$ coincides with the Dirichlet distribution $\BD_\theta$. Recall that    $ \BG_\theta$   is the unique probability measure on $\MM$ such that for any finitely many disjoint measurable subsets $\{A_1,\cdots, A_n\}$ of $M$,
$$\MM\ni\eta\mapsto \eta(A_i),\ \ 1\le i\le n$$ are independent Gamma random variables with shape parameters $\{\theta(A_i)\}_{1\le i\le n}$ and scale parameter $1$; that is,
\beq\label{*0}  \int_{\MM} f(\eta(A_1),\cdots, \eta(A_n))\BG_\theta(\d\eta)
 = \int_{[0,\infty)^n} f(s_1,\cdots, s_n)\prod_{i=1}^n \gg_{\theta(A_i)}(\d s_i) \end{equation}  holds for any $ f\in \B_b(\R^n),$  where for a constant $r>0$,
 \beq\label{*0'}  \gg_r(\d s):= 1_{[0,\infty)}(s) \ff{s^{r-1}\e^{-s}}{\GG(r)} \, \d s,\ \  \GG(r):=\int_0^\infty s^{r-1}\e^{-s}\d s,\end{equation} and we set $\gg_0= \dd_0$,   the Dirac measure at point $0$.

 The Gamma distribution   $\BG_\theta$ is supported  on the class of  finite discrete measures
 $$\MM_{dis}:=\Big\{\sum_{i=1}^\infty s_i\dd_{x_i}:\ s_i\ge 0, x_i\in M, \sum_{i=1}^\infty s_i<\infty\Big\}.$$
Moreover,  under $\BG_\theta$ the random variables $\eta(M)\in (0,\infty)$ and $\Psi(\eta):= \ff{\eta}{\eta(M)}\in \MM_1$ are independent with
\beq\label{BDD0}\BG_\theta\big(\eta(M)<r, \Psi(\eta)\in \mathbf A)=    \ff {\mathbb D_\theta(\mathbf A)} {\GG(\theta(M))} \int_0^r s^{\theta(M)-1}\e^{-s}\d s,\ \ r>0, \mathbf A\in \B(\MM_1).\end{equation}
Consequently,
\beq\label{BDD} \mathbb D_\theta= \BG_\theta\circ \Psi^{-1},\ \ \Psi(\eta):= \ff{\eta}{\eta(M)},\ \eta\in\MM\setminus\{0\}.\end{equation}

Both $\BD_\theta$ and $\BG_\theta$ are images of the Poisson measure $\pi_{\hat\theta}$ with intensity
\beq\label{HTT}  \hat\theta(\d x,\d s):= s^{-1}\e^{-s}\theta(\d x)\d s \end{equation}  on the product manifold $\hat M:= M\times (0,\infty)$. Recall that $\pi_{\hat\theta}$  is the unique probability measure on the configuration space
$$ \GG_{\hat M}:=\Big\{\gg=\sum_{i=1}^\infty \dd_{(x_i,s_i)}:\ (x_i,s_i)\in \hat M,   \gg(K)<\infty\ \text{for\ any\ compact}\
K\subset \hat M\Big\} $$ equipped with the vague topology, such that for any disjoint compact subsets $\{K_i\}_{1\le i\le n}$ of $\hat M$, $\gg\mapsto \gg(K_i)$ are independent Poisson random variables of parameters $\hat\theta(K_i)_{1\le i\le n}.$   By    \cite[Theorem 6.2]{HKPR}, we have
\beq\label{GGT0}\s(\gg):= \sum_{i=1}^\infty s_i<\infty,\ \ \text{for}\ \pi_{\hat\theta}\text{-a.s.}\ \gg=\sum_{i=1}^\infty\dd_{(x_i,s_i)}\in \GG_{\hat M},\end{equation}
 and
\beq\label{GGT} \BG_\theta  = \pi_{\hat\theta}\circ\Phi^{-1},\end{equation}   where
$$\Phi(\gg):= \sum_{i=1}^\infty s_i \dd_{x_i}\in \MM,\ \ \gg:=\sum_{i=1}^\infty \dd_{(x_i,s_i)}\in \GG_{\hat M} \ \text{with}\ \sum_{i=1}^\infty s_i<\infty.$$
Combining \eqref{BDD} with \eqref{GGT}, we obtain
\beq\label{BDD2} \mathbb D_\theta= \BG_\theta\circ\Psi^{-1} =\pi_{\hat\theta}(\Psi\circ\Phi)^{-1}.\end{equation}

\subsection{Intrinsic and extrinsic derivatives  }

These derivatives were introduced in \cite{AKR} and \cite{ORS} on the configuration space and the space of probability measures respectively, which can be extended to  $\MM$ under the map $\Phi: \GG_{\hat M}\to \MM$, see for instance \cite{KLV}.

To introduce the intrinsic derivative for a function on $\MM$ (or $\MM_1$), we let $\scr V_0(M)$ be the class of smooth vector fields  with compact supports on $M$. For any $v\in \scr V_0(M)$, let
$$\phi_v(x)= \exp_x[v(x)],\ \ x\in M,$$ where $\exp$ is the exponential map on $M$. Then $\phi_v\in C^\infty(M\to M)$. For a function $F$ on $\MM$, we define its directional derivative  along $v$ by
$$\nn^{int}_vF(\eta):=\lim_{\vv\downarrow 0} \ff{F(\eta\circ \phi_{\vv v}^{-1})- F(\eta)}\vv$$
if it exists. Let $L^2(\scr V(M),\eta)$ be the space of   all measurable vector fields $v$ on $M$ with $\eta(|v|^2)<\infty$.  When $\nn^{int}_v F(\eta)$ exists for all $v\in \scr V_0(M)$ such that
$$|\nn^{int}_v F(\eta)|\le c \|v\|_{L^2(\eta)},\ \ v\in \scr V_0(M)$$ holds for some constant $c\in (0,\infty)$, then by Riesz representation theorem there exists a unique $\nn^{int}F(\eta)\in L^2(\scr V(M),\eta)$   such that
\beq\label{DD0} \nn^{int}_vF(\eta)= \<\nn^{int}F(\eta), v\>_{L^2(\eta)}:=\int_M  \<\nn^{int}F(\eta), v\>_M\,\d\eta,\ \ v\in \scr V_0(M).\end{equation} In this case, we call $F$  intrinsically  differentiable  at $\eta$ with derivative $\nn^{int}F(\eta)$. If $F$ is intrinsically differentiable at all $\eta\in \MM$ (or $\MM_1$), we call it intrinsically differentiable on $\MM$ (or $\MM_1$).

Next, a measurable real function $F$ on $\MM$   is called extrinsically differentiable at $\eta\in \MM$,    if
 $$\nn^{ext} F(\eta)(x):= \ff{\d}{\d s} F(\eta+s\dd_x)\Big|_{s=0}\ \text{exists\ for\ all} \   x\in M,$$ such that
 $$\|\nn^{ext} F(\eta)\|:= \|\nn^{ext} F(\eta)(\cdot)\|_{L^2(\eta)} <\infty.$$   When a function $F$ on $\MM_1$ is considered, it is called intrinsically differentiable if 
   $$\tt\nn^{ext} F(\mu)(x):= \ff{\d}{\d s} F((1-s)\mu+s\dd_x)\Big|_{s=0}\ \text{exists\ for\ all} \   x\in M,\mu\in \MM_1$$ and $\tt\nn^{ext} F(\mu)\in L^2(\mu)$. 
   If $F$ is   extrinsically differentiable at all $\eta\in \MM$ (or $\MM_1$),  we    call it extrinsically differentiable on $\MM$ (or $\MM_1$).
Let $\D(\MM)$ (respectively  $\D(\MM_1)$) denote the set of functions which are  intrinsically  and extrinsically differentiable on $\MM$ (respectively $\MM_1$).

A typical subclass of $\D(\MM)$ and $\D(\MM_1)$  is the set of cylindrical functions
\beq\label{CLL} \F C_b^\infty(M):= \big\{\eta\mapsto f(\<h_1,\eta\>,\cdots, \<h_n,\eta\>):\ n\ge 1, f\in C_b^\infty(\R^n), h_i\in C_0^\infty(M)\big\},\end{equation}  where $\<h_i,\eta\>:=\eta(h_i)=\int_M h_i\d\eta.$
This class is dense in $L^2(\MM_1,\mathbb D_\theta)$ and $L^2(\MM, \BG_\theta)$, and the cylindrical function $F:=f(\<h_1,\cdot\>,\cdots, \<h_n,\cdot\>)$ is differentiable with
\beq\label{DD} \beg{split} &\nn^{int} F(\eta)= \sum_{i=1}^n (\pp_i f)(\<h_1,\eta\>,\cdots, \<h_n,\eta\>) \nn h_i,\\
  &\nn^{ext} F(\eta)=  \sum_{i=1}^n (\pp_i f)(\<h_1,\eta\>,\cdots, \<h_n,\eta\>)   h_i, \ \ \eta\in\MM, \end{split} \end{equation} where $\nn$ is the gradient operator on $M$.
Restricting on $\MM_1$ we will consider
$$\tt\nn^{ext} F(\mu):=  \sum_{i=1}^n (\pp_i f)(\<h_1,\mu\>,\cdots, \<h_n,\mu\>)   (h_i-\mu(h_i)), \ \ \mu\in\MM_1,$$
which is the centered extrinsic derivative of $F$ at $\mu$  since
\beq\label{DDW} \mu\big(\tt\nn^{ext}F(\mu)\big)=0,\ \ \mu\in \MM_1.\end{equation}
See \cite{RW19} for general results on the relations of $\nn^{int},\nn^{ext}$ and $\tt\nn^{ext}.$ 
\subsection{The main result}

Now, for any $\ll>0$,  we consider the square fields for $F,G\in \F C_b^\infty(M)$:
\beq\label{SQF} \beg{split} &\GG^\ll(F,G)(\eta) := \int_M \Big\{\<\nn^{int}F(\eta), \nn^{int}G(\eta)\>_M + \ll   (\nn^{ext}F(\eta))   \nn^{ext}G(\eta) \Big\}(x)\,\eta(\d x),\ \ \eta\in\MM, \\
&\tt\GG^\ll(F,G)(\mu) := \int_M \Big\{\<\nn^{int}F(\mu), \nn^{int}G(\mu)\>_M + \ll  (\tt\nn^{ext}F(\mu))  \tt \nn^{ext}G(\mu) \Big\}(x)\,\mu(\d x),\ \ \mu\in\MM_1,\end{split}\end{equation}
which lead to   the following bilinear forms on $L^2(\MM_1,\BD_\theta)$ and $L^2(\MM,\BG_\theta)$ respectively:
\beq\label{EENU} \beg{split} &\EE_{\BD_\theta}^\ll (F,G):=\int_{\MM_1} \tt\GG^\ll(F,G)(\mu)\, \mathbb D_\theta(\d\mu),\\
&\EE_{\BG_\theta}^\ll (F,G):=\int_{\MM} \GG^\ll(F,G)(\eta)\, \BG_\theta(\d\eta),\ \ F,G\in \F C_b^\infty(M).\end{split} \end{equation}
To ensure the closability of these bilinear forms, we take
\beq\label{THETA} \theta(\d x)= \e^{V(x)} {\rm vol}(\d x)\ \text{for\ some}\ V\in W_{loc}^{1,1}(M),\ \ \theta(M)<\infty, \end{equation} where ${\rm vol}$ is the Riemannian  volume measure. Then the integration by parts formula gives
\beq\label{SY1} \EE_\theta(h_1,h_2):= \int_M \<\nn h_1,\nn h_2\>_M(x)\,\theta(\d x)=-\int_M h_1(\DD+\nn V)h_2 \d\theta, \ \ h_1,h_2\in C_0^\infty(M).\end{equation}
So, the bilinear form is closable in $L^2(M,\theta)$, and the closure $(\EE_\theta,\D(\EE_\theta))$ is a Dirichlet form.

We will prove the closability of $(\EE_{\BG_\theta}^\ll,\F C_b^\infty(M))$ and $(\EE_{\BD_\theta}^\ll,\F C_b^\infty(M))$, and  calculate the   spectral gaps for the corresponding Dirichlet forms.

Recall that for  a probability space $(E,\scr B,\mu)$ and a symmetric Dirichlet form $(\EE,\D(\EE))$ on $L^2(E,\mu)$   with $1\in \D (\EE)$ and $\EE(1,1)=0$, the spectral gap  of the Dirichlet form  is given by
$$\gap(\EE)= \inf\Big\{\EE (F,F):\ F\in \D(\EE), \mu(F)=0, \mu(F^2)=1\Big\}.$$
 By the spectral theorem,
$\gap(\EE)$ is the exponential convergence rate of the associated Markov semigroup $(P_t)_{t\ge 0}$, i.e.
$$\|  P_t-\mu\|_{L^2(\mu)} := \sup_{\mu(F^2)\le 1} \|P_t F -\mu(F)\|_{L^2(\mu)}=\e^{-\gap(\EE)t},\ \ t\ge 0.$$

 Let
 $$\ll_\theta =\gap(\EE_\theta):=\inf\big\{\theta(|\nn f|^2):\ f\in C_b^1(M), \theta(f)=0, \theta(f^2)=1\big\}$$ be the spectral gap of the Dirichlet form $(\EE_\theta,\D(\EE_\theta))$ in $L^2(M,\theta)$.
The main result of this paper is the following.

\beg{thm}\label{T1.1} Let $\theta$ be in $\eqref{THETA}$ and let $F_0(\eta)=\eta(M), \eta\in \MM$.  \beg{enumerate}
\item[$(1)$] $(\EE_{\BG_\theta}^\ll, \F C_b^\infty(M))$ is closable in $L^2(\MM,\BG_\theta)$ and the closure $(\EE_{\BG_\theta}^\ll, \D(\EE_{\BG_\theta}^\ll)) $ is a quasi-regular symmetric Dirichlet form with spectral gap $\gap(\EE^\ll_{\BG_\theta})=\ll.$
\item[$(2)$]  $(\EE_{\BD_\theta}^\ll, \F C_b^\infty(M))$ is closable in $L^2(\MM_1,\BD_\theta)$ and the closure $(\EE_{\BD_\theta}^\ll, \D(\EE_{\BD_\theta}^\ll)) $ is a quasi-regular  symmetric Dirichlet form with spectral gap satisfying
    $$ 
    \ll \theta(M)\le   \gap(\EE^\ll_{\BD_\theta})\le \ll\theta(M)+\ll_\theta(\theta(M)+1).$$   Consequently, if   $\ll_\theta=0,$ then $\gap(\EE^\ll_{\BD_\theta})= \ll\theta(M).$ Moreover, $F_0(F\circ\Psi)\in \D(\EE_{\BG_\theta}^\ll)$ for $F\in \F C_b^\infty(M)$ and for any $F,G\in \F C_b^\infty(M)$,
    \beq\label{*DR}  \EE_{\BG_\theta}^\ll(F_0(F\circ\Psi),F_0(G\circ\Psi))= \theta(M) \EE_{\BD_\theta}^\ll(F,G)+\ll\theta(M)\BD_\theta(F,G).\end{equation} \end{enumerate} \end{thm}

 The formula \eqref{*DR} will play a crucial role in the proof of the closability of $(\EE_{\BD_\theta}^\ll,\F C_b^\infty(M))$.
When $\ll_\theta>0$, the exact value of $\gap(\EE^\ll_{\BD_\theta})$ is unknown. Since in this case the intrinsic derivative part will play a non-trivial role,   we believe that
$ \gap(\EE^\ll_{\BD_\theta})$ is strictly larger than $\ll\theta(M),$ hopefully our upper bound could be sharp.

When $\ll=1$ and without the intrinsic derivative part, the Dirichlet form $\EE_{\BD_\theta}^\ll$ reduces to
\beq\label{FV} \EE_{\BD_\theta}^{FV} (F,G):=   \int_{\MM_1}\BD_\theta(\d\mu) \int_M   \big\{(\tt\nn^{ext}F(\mu))   \tt\nn^{ext}G(\mu) \big\}\d\mu,\end{equation} which is associated with the   Fleming-Viot process.
 It  has been derived in \cite{SN} that  \beq\label{FVG} \gap(\EE_{\BD_\theta}^{FV} )=\theta(M),\end{equation}  see   \cite{STA} for the study of log-Sobolev inequality in finite-dimensions as well as \cite{FMW,W18,WZ18} for  functional inequalities of a modified Fleming-Viot process.

When $\ll=1$ and $V=0$,  \cite[Theorem 14]{KLV} presents an integration by parts formula   for $\EE_{\BG_\theta}^\ll$ in the class
$$\tt\F C_b^{\infty} (M):= \big\{F\circ \psi_\vv:\ F\in \F C_b^\infty(M),\vv>0\big\},$$
where
$$\psi_\vv(\eta)= \sum_{\eta(\{x\})\ge \vv} \eta(\{x\})\dd_x,\ \ \vv>0, \eta\in\MM.$$
Therefore, letting $(\L_{\BG_\theta}^\ll, \D(\L_{\BG_\theta}^\ll))$ be the generator of the Dirichlet form $(\EE_{\BG_\theta}^\ll, \D(\EE_{\BG_\theta}^\ll))$, we have $\D(\L_{\BG_\theta}^\ll)\supset \tt\F C_b^\infty(M).$ However, in general $\F C_b^\infty(M)$ is not included in $\D(\L_{\BG_\theta}^\ll)$.
 Indeed, according to \cite{KLV} we have the integration by parts formula
 \beq\label{*DP}\int_\MM \nn^{int}_v(F\circ\psi_\vv)\, \d\BG_\theta =- \int_\MM (F\circ\psi_\vv)  B_{\vv,v}  \, \BG_\theta(\d\eta),\ \ v\in \scr V_0(M), F\in \F C_b^\infty(M) \end{equation}  for  $B_{\vv,v}(\eta):= \sum_{\eta(\{x\})\ge \vv}  \big\{\div (v)+\<v,\nn V\>\big\}(x)$, where $\div$ is the divergence operator in $M$. This formula makes sense because
 $$| \BG_\theta(B_{\vv,v})|  =  |\theta( \div+\<v,\nn V(x)\>_M)|  \int_\vv^\infty s^{-1}\e^{-s}\d s<\infty.$$
 So,   $F\circ \psi_\vv\in \D(\L_{\BG_\theta}^\ll)$ for any $\vv>0$ and $F\in \F C_b^\infty(M)$.  However, since $ \int_0^\infty s^{-1}\e^{-s}\d s=\infty$, \eqref{*DP} does not make sense for $\vv=0$.

\

To prove Theorem \ref{T1.1}, we will formulate the bilinear form $\EE_{\BD_\theta}$ as the image of the Dirichlet form on the configuration space constructed in \cite{AKR}, for which the (weak) Poincar\'e inequality has been established in \cite[Section 7]{RW01}. To this end, we first recall in Section 2 some known  results on the configuration space, then prove Theorem \ref{T1.1}(1) in Section 3 by transforming these results to  the Gamma process on $\MM$, and finally prove Theorem \ref{T1.1}(2) in Section 4 by mapping the Gamma process to the subclass $\MM_1$.

\

\section{Analysis on the configuration space}

In this section, we first recall the diffusion process on the configuration space constructed in \cite{AKR,AKR2}, then calculate  the spectral gap.

For $F:=f(\<\hat h_1,\cdot\>,\cdots, \<\hat h_n, \cdot\>)\in \F C_b^\infty(\hat M),$ where $f\in C_b^\infty(\R^n)$ and $\hat h_i\in C_0^\infty(\hat M)$, let
\beq\label{DD'}\nn^\GG F(\gg)= \sum_{i=1}^n (\pp_if)(\<\hat h_1,\gg\>,\cdots, \<\hat h_n,\gg\>) \hat \nn \hat h_i,\ \ \gg\in \GG_{\hat M}.\end{equation} For $\ll>0$, we take
  the following Riemannian metric on the manifold $\hat M:= M\times (0,\infty):$
\beq\label{*OP} \<a_1\pp_s +v_1, a_2\pp_s + v_2\>_{\hat M} := (\ll s)^{-1} a_1a_2 + s \<v_1,v_2\>_M,\ \ a_1,a_2\in\R, v_1, v_2\in TM.\end{equation}
Let $\hat \DD,\hat \nn$ and $\hat{\rm vol}$ be the Laplacian, gradient and volume measure on $\hat M$ respectively.

 Consider the bilinear form
 \beq\label{GMMP} \EE_{\hat\theta}^\GG(F,G) := \int_{\GG_{\hat M}}   \pi_{\hat\theta}(\d\gg) \int_{\hat M} \<\nn^\GG (F\circ\Phi)(\gg), \nn^\GG (G\circ\Phi)(\gg)\>_{\hat M}\d\gg,\ \ F,G\in \F C_b^\infty(M).\end{equation}
To formulate the integration by parts formula of this form,
  we intend to find out a function $W$ on $\hat M$ such that
\beq\label{WW} \e^{W(s,x)} \hat{\rm vol}(\d s,\d x)= \hat\theta(\d s,\d x), \end{equation} where $\hat\theta(\d s,\d x):= s^{-1}\e^{-s} \d s\theta(\d x)$ is given by \eqref{HTT}. So, for
\beq\label{HATL} \hat L f:= \hat\DD f +\<\hat\nn W, \hat\nn f\>_{\hat M},\ \ f\in C^2(\hat M),\end{equation}  we have the integration by parts formula
\beq\label{ITL} \EE_{\hat  \theta}(\hat h_1,\hat h_2):= \int_{\hat M} \<\hat \nn\hat h_1, \hat\nn\hat h_2\>_{\hat M} \d\hat\theta=- \int_{\hat M} (\hat h_1\hat L\hat h_2)\d\hat\theta,\ \ \
\hat h_1,\hat h_2\in C_0^\infty(\hat M).  \end{equation}
Therefore, letting
 \beq\label{GTGG}\beg{split} \L_{\hat \theta}^\GG F(\gg)= &\sum_{i,j=1}^n (\pp_i \pp_j f)(\<\hat h_1,\gg\>,\cdots,\<\hat h_n,\gg\>) \gg(\<\hat \nn\hat h_i,\hat\nn\hat h_j\>_{\hat M})\\
&+ \sum_{i=1}^n (\pp_i  f)(\<\hat h_1,\gg\>,\cdots,\<\hat h_n,\gg\>) \gg(\hat L\hat h_i),\ \ \gg\in \GG_{\hat M},\end{split}\end{equation}
we have
 (see  \cite[Theorem 4.3]{AKR2})
\beq\label{TTD}   \EE_{\hat\theta}^\GG(F,G)
  =-\int_{\GG_{\hat M}} (G\L_{\hat\theta}^\GG F)\, \d \pi_{\hat\theta},\ \ F,G\in \F C_b^\infty(\hat M),\end{equation} which implies the closability of $( \EE_{\hat\theta}^\GG,\F C_b^\infty(\hat M))$, so that the closure $( \EE_{\hat\theta}^\GG,\D(\EE_{\hat\theta}^\GG))$ is a symmetric Dirichlet form in $L^2(\GG_{\hat M},\pi_{\hat\theta})$. Moreover, as explained in \cite[Section 4.5.1]{MR},   the  result \cite[Corollary 4.9]{MR} applies to this situation, so that the Dirichlet form $( \EE_{\hat\theta}^\GG,\D(\EE_{\hat\theta}^\GG))$ is quasi-regular and local, and hence is associated with a diffusion process on $\GG_{\hat M}$. Recall that $\GG_{\hat M}$ is equipped with the vague topology.

To calculate the generator $\L_{\hat\theta}^\GG$ defined in \eqref{GTGG}, it suffices to figure out the operator $\hat L.$ To this end,
for a fixed point $z:=(\bar s, \bar x)\in \hat M$,  we take the normal coordinates $(s,x_1,\cdots, x_t)$ in a neighbourhood $\scr O(z)$ of $z$ such that
$$  U_1:=\pp_s, \ \ U_{i+1}:=\pp_{x_i},\ 1\le i\le d $$ satisfy
\beq\label{**7}  g(z) = {\rm diag}\{(\ll \bar s)^{-1}, \bar s,\cdots, \bar s\},\ \
  \hat\nn U_i(z)=0,\ 1\le i\le d+1.\end{equation}  Let
 $$ g(s,x):= (\<U_i,U_j\>_{\hat M})_{1\le i,j\le d+1}(s,x),\ \ (s,x) \in\scr O(z).$$
Then the Riemannian volume measure on $\scr O(z)$ is
$$\hat{\rm vol} (\d s,\d x)= \ss{{\rm det} g(s,x)}\, \d s \d x,$$ so that  \eqref{WW} holds for
$$W(s,x):=\log\Big[\ff{\hat\theta(\d s,\d x)}{\hat{\rm vol}(\d s,\d x)}\Big]= -s-\log s +V(x) -\ff 1 2   \log \big[{\rm det} g(s,x)\big],\ \ (s,x)\in \scr O(z).$$ Letting $(g^{ij})_{1\le i,j\le d+1}=g^{-1}$, we derive from \eqref{**7} that $\nn g^{ij}(z)=0$. So,
\beq\label{**9} \hat\nn W(z):=\sum_{i=1}^{d+1} \big\{g^{ii} (U_i W) U_i\big\}(z) = \bar s^{-1} \nn V(\bar x)-\ll(1+\bar s) \pp_s.\end{equation} By the same reason,
  for any $\hat h\in C^\infty(\hat M)$ we have
\beg{align*} \hat\nn \hat h (z)&:= \sum_{i=1}^{d+1} \big\{g^{ii} (U_i\hat h) U_i\big\}(z) = \ll \bar s (\pp_s \hat h)(\bar s,\bar x)\pp_s + \bar s^{-1} \sum_{i=1}^d( \pp_{x_i}\hat h(\bar s,\bar x) )\pp_{x_i},\\
 \hat \DD \hat h(z)&:= \ff 1 {\ss{{\rm det} g(z)}} \sum_{i=1}^{d+1} U_i\Big(\ss{{\rm det}g} g^{ii} U_i \hat h\Big)(z)= \sum_{i=1}^{d+1} \big\{g^{ii} U_i^2 \hat h\big\}(z)\\
&  = \ll \bar s \pp_s^2 \hat h(\bar s,\bar x)+ \bar s^{-1}\DD\hat h(\bar s,\cdot)(\bar x).\end{align*}
This together with \eqref{**9} implies that at point $z$,
\beq\label{*PY} \beg{split} \hat L&:= \hat\DD +\hat \nn W=\ll s(\pp_s^2 -\pp_s) -\ll\pp_s+s^{-1}(\DD+ \nn V).\end{split} \end{equation} Since $z\in\hat M$ is arbitrary, this formula holds for all points $(s,x)\in \hat M.$



\beg{thm} \label{T2.1}  Let $\ll>0$ and $ \theta\in\MM$ be as in $\eqref{THETA}$. Then $\gap(\EE_{\hat\theta}^\GG)=\ll.$    \end{thm}

\beg{proof} According to \cite{SIM}, see also \cite[Theorem 7.1]{RW01}, we have
\beq\label{SPL} \gap (\EE_{\hat\theta}^\GG)= \ll_{\hat\theta}:= \inf\big\{\EE_{\hat\theta}(\hat h,\hat h): \hat h\in \D(\EE_{\hat \theta}), \hat\theta(\hat h^2)=1\big\},\end{equation}
where $ \D(\EE_{\hat \theta})$ is the closure of $C_0^\infty(\hat M)$ under the Sobolev norm $\|\hat h\|_{1}:= \ss{\hat\theta(|\hat h|^2)+ \EE_{\hat\theta}(\hat h,\hat h)}.$ Let
 $$\hat h(x,s):= s+1,\ \ (x,s)\in\hat M.$$  By \eqref{*PY} we have
\beq\label{*PY2}\hat L \hat h(x,s)= -\ll \hat h (x,s),\ \ (x,s)\in\hat M.\end{equation}
Combining this with \eqref{ITL},   for any $g\in C_0^\infty(\hat M)$ we have
\beg{align*} & \ll\hat\theta(g^2)= -\int_{\hat M} \ff {g^2} {\hat h} \hat L \hat h\,\d\hat\theta =\int_{\hat M} \big<\hat\nn (g^2/\hat h),\hat\nn \hat h\big\>_{\hat M} \d\hat\theta\\
& = \int_{\hat M}\big\{2g\<\hat\nn g,\hat\nn \log\hat h\>_{\hat M} -g^2\<\hat\nn\log \hat h, \hat\nn \log \hat h\>_{\hat M} \big\}\d\hat\theta\\
 &\le \int_{\hat M}  \<\hat\nn g, \hat\nn g\>_{\hat M}  \d\hat\theta=\EE_{\hat\theta}(g,g).\end{align*}
Therefore,
$  \ll_{\hat\theta}\ge \ll.$

On the other hand,
since $M$ is complete, there exists a sequence $\{h_n\}_{n\ge 1}\subset C_0^\infty(M)$ such that
\beq\label{HNN} 0\le h_n\uparrow 1\ \text{as}\ n\uparrow\infty,\ \text{and}\ \|\nn h_n\|_\infty\le \ff 1 n,\ \ n\ge 1.\end{equation}
For any $\vv\in (0,1)$ let
$$\hat h_{n,\vv}(x,s)= (s-\vv)^+ h_n(x),\ \ \hat h_\vv(x,s)= (s-\vv)^+,\ \ (x,s)\in \hat M, n\ge 1.$$
Then $\{\hat h_{n,\vv}\}_{n\ge 1,\vv\in (0,1)}\subset \D(\EE_{\hat\theta})$ with $0\le \hat h_{n, \vv}\le \hat h_{\vv}.$ So, for fixed $\vv$, by the dominated convergence theorem
$$\lim_{n\to\infty} \int_{\hat M} |\hat h_{n,\vv}-\hat h_\vv|^2\d\hat\theta =\int_{\hat M} \Big(\lim_{n\to\infty}   |\hat h_{n,\vv}-\hat h_\vv|^2\Big)\d\hat\theta =0,$$ and due to \eqref{ITL} and \eqref{HNN},
\beg{align*} &\lim_{n,m\to\infty} \EE_{\hat\theta}(\hat h_{n,\vv}-\hat h_{m,\vv}, \hat h_{n,\vv}-\hat h_{m,\vv})\\
 &= \lim_{n,m\to\infty} \int_{\hat M} \big\{\ll s 1_{\{s\ge \vv\}} |h_n-h_m|^2(x)+ s^{-1} |(s-\vv)^+|^2 |\nn(h_n-h_m)|^2(x)\big\}s^{-1}\e^{-s}\d s\theta(\d x)\\
&\le  \lim_{n,m\to\infty}   \Big(\ll \theta(|h_n-h_m|^2) + \ff{\theta(M)}{n^2\land m^2}\Big) =0.\end{align*}
Thus, $\hat h_\vv\in \D(\EE_{\hat\theta})$ with
\beq\label{RPP} \EE_{\hat\theta}(\hat h_\vv,\hat h_\vv)= \lim_{n\to\infty} \EE_{\hat\theta}(\hat h_{n,\vv},\hat h_{n,\vv}) =\ll\theta(M)\int_\vv^\infty  \e^{-s}\d s,\ \ \vv\in (0,1).\end{equation}
Combining this with
$$\lim_{\vv\to 0} \hat\theta(h_\vv^2)=\theta(M) \int_0^\infty s \e^{-s}\d s=\theta(M),$$
we obtain
$$\ll_{\hat \theta}\le \lim_{\vv\downarrow 0} \ff{\EE_{\hat\theta}(\hat h_\vv,\hat h_\vv)}{\hat\theta(h_\vv^2)}=\ll.$$
This together with $\ll_{\hat \theta}\ge \ll$ derived above gives $\ll_{\hat \theta}=\ll.$ So, the proof is finished by \eqref{SPL}.   \end{proof}

 \section{Proof of Theorem \ref{T1.1}(1)}

\beg{thm}\label{T3.1}  Let $\theta$ be as in $\eqref{THETA}$. Then
 $(\EE_{\BG_\theta}^\ll, \F C_b^\infty (M))$ is closable in $L^2(\MM,\BG_\theta)$ and the closure $(\EE_{\BG_\theta}^\ll, \D(\EE_{\BG_\theta}^\ll))$ is a quasi-regular  symmetric Dirichlet form.
Moreover,
$ \gap(\EE_{\BG_\theta})=\ll. $ \end{thm}

\beg{proof} Due to \eqref{GGT}, we may prove this result by using \eqref{TTD} and Theorem \ref{T2.1}.
For any $h\in C_0^\infty(M)$, let $\hat h(x,s)=sh(x), (x,s)\in \hat M.$ Then
 \beq\label{PHIG} \Phi(\gg)(h )= \gg(\hat h),\ \ \gg\in \GG_{\hat M}.\end{equation}
Let $F=f(\<h_1,\cdot\>,\cdots,\<h_n,\cdot\>), G= g(\<h_1,\cdot\>,\cdots,\<h_n,\cdot\>)\in \F C_b^\infty(M)$. We have
$$F\circ \Phi= f(\<\hat h_1,\cdot\>,\cdots, \<\hat h_n,\cdot\>), \ G\circ \Phi= g(\<\hat h_1,\cdot\>,\cdots, \<\hat h_n,\cdot\>)\in\D(\EE_{\hat\theta}^\GG)$$
such that
 $$\<\nn^\GG (F\circ\Phi), \nn^\GG (G\circ\Phi)\>_{\hat M}(\gg)  = \sum_{i,j=1}^n \big\{(\pp_i f)(\pp_j g)\big\}(\<\hat h_1, \gg\>,\cdots, \<\hat h_n, \gg\>)\gg\big(\<\hat\nn \hat h_i, \hat\nn \hat h_j\>_{\hat M}\big).$$ Because of \eqref{*OP} and  \eqref{PHIG}, we have
 \beq\label{*PYN}\beg{split}  \gg\big(\<\hat\nn \hat h_i, \hat\nn \hat h_j\>_{\hat M}\big)& = \int_{\hat M} \big\{\ll s (h_ih_j)(x) + s \<\nn h_i,\nn h_j\>_M(x)\big\}\gg(\d x,\d s)\\
 &= \Phi(\gg) \big(\<\nn   h_i,  \nn   h_j\>_{ M}+\ll h_ih_j\big).\end{split} \end{equation}
 Thus,
 \beg{align*}&\int_{\hat M} \<\nn^\GG (F\circ\Phi)(\gg), \nn^\GG (G\circ\Phi)(\gg)\>_{\hat M}\d\gg \\
 &=  \sum_{i,j=1}^n \big\{(\pp_i f)(\pp_j g)\big\}(\<h_1, \Phi(\gg)\>,\cdots, \<h_n, \Phi(\gg)\>)\Phi(\gg)\big(\<\nn   h_i,  \nn   h_j\>_{ M}+\ll h_ih_j\big)\\
 &= \GG^\ll(F,G)(\Phi(\gg)),\ \ \gg\in \GG_{\hat M}.\end{align*}
Combining this with  \eqref{GGT}, \eqref{EENU} and \eqref{GMMP}, we obtain
\beq\label{*PY3}  \EE_{\hat\theta}^\GG(F\circ \Phi, G\circ\Phi)   = \int_{\MM}  \GG^\ll(F,G) \, \d\BG_\theta = \EE_{\BG_\theta}^\ll(F,G).  \end{equation}
Below  we   prove the closability, quasi-regularity, and   the spectral gap bounds respectively.

(a) The closability. Let $\{F_n\}_{n\ge 1}\subset \F C_b^\infty(M)$ such that
\beq\label{NMM} \lim_{n\to\infty} \BG_\theta(F_n^2)= 0\ \text{and}\ \lim_{n,m\to\infty}\EE_{\BG_\theta}^\ll(F_n-F_m, F_n-F_m) =0.\end{equation}
It remains to show that $\lim_{n\to\infty} \EE_{\BG_\theta}^\ll(F_n,F_n)=0.$
By \eqref{GGT} and \eqref{*PY3}, \eqref{NMM} implies
$$\lim_{n\to\infty} \pi_{\hat \theta}(|F_n\circ\Phi|^2)+ \lim_{n,m\to\infty}\EE_{\hat\theta}^\GG(F_n\circ\Phi-F_m\circ\Phi, F_n\circ\Phi-F_m\circ\Phi) =0,$$
so that the closability of $(\EE_{\hat\theta}^\GG,\D(\EE_{\hat\theta}^\GG))$ and \eqref{NMM} imply
$$ \lim_{n\to\infty} \EE_{\BG_\theta}^\ll(F_n,F_n)= \lim_{n \to\infty}\EE_{\hat\theta}^\GG(F_n\circ\Phi, F_n\circ\Phi) =0.
$$

(b) The quasi-regularity. According to \cite[Chap. IV, Def. 3.1]{MRB}, the Dirichlet  form $(\EE_{\BG_\theta}^\ll, \D(\EE_{\BG_\theta}^\ll))$  is quasi-regular if and only if there exist a sequences of compact subsets $\{K_n\}_{n\ge 1}$ of $\MM$ such that the class
\beq\label{KKL} \D_K:= \big\{F\in \D(\EE_{\BG_\theta}^\ll): F|_{\MM\setminus K_n}=0\ \text{for\ some\ }n\ge 1\big\}\end{equation}
is dense in $\D(\EE_{\BG_\theta}^\ll)$ under the Sobolev norm
$$ \|F\|_{L^2(\BG_\theta)}+ \ss{\EE_{\BG_\theta}^\ll(F,F)},\ \ F\in \D(\EE_{\BG_\theta}^\ll).$$
In this case, the sequence $\{K_n\}$ is called a $\EE_{\BG_\theta}^\ll$-nest.

To construct such a $\EE_{\BG_\theta}^\ll$-nest, we choose a function $\psi\in C^\infty([0,\infty))$ with $1\le \psi(r)\uparrow\infty$ as $r\uparrow\infty$ and $0\le \psi'\le 1,$ such that $\theta(|\psi\circ\rr_o|^2)<\infty$, where $\rr_o$ is the Riemannian distance from a fixed point $o\in M$.
 Thus,
\beq\label{PPF} f:= \psi\circ \rr_o\in \D(\EE_\theta).\end{equation}  Since $M$ is complete and $\psi\circ\rr_o\uparrow\infty$ as $\rr_o\uparrow\infty$, the level sets
 $$K_n:=\{\eta\in\MM: \eta(\psi\circ\rr_o)\le n\},\ \ n\ge 1$$ are compact in $\MM$.  It remains to show that $\{K_n\}_{n\ge 1}$ is a $\EE_{\BG_\theta}^\ll$-nest.
Since $\F C_b^\infty(M)$ is dense in $\D(\EE_{\BG_\theta}^\ll)$, it suffices to show that for any $F\in \F C_b^\infty(M)$, there exist a sequence
$\{F_n\}_{n\ge 1}\subset \D_K$, where $\D_K$ is defined in \eqref{KKL} for the present $\{K_n\}_{n\ge 1}$,  such that
\beq\label{NNN} \lim_{n\to\infty}    \Big\{\BG_\theta(|F_n-F|^2)+\EE_{\BG_\theta}^\ll(F_n-F,F_n-F)\Big\} =0.\end{equation}
We will prove this formula for
\beq\label{LLFN} F_n:= F\cdot \big\{1 \land (n+1-F_0)^+\big\},\ \ n\ge 1, F_0(\eta):= \eta(\psi\circ\rr_o).\end{equation}
To this end, we first confirm  that $ F_0\in \D(\EE_{\BG_\theta}^\ll)$, so that
$\{F_n\}_{n\ge 1}\subset \D_K$ by the definitions of $K_n$ and $F_n$.
By \eqref{PPF},  there exist functions  $\{f_n\}_{n\ge 1}\subset C_0^\infty(M)$ such that
\beq\label{*09}0\le f_n\le f:= \psi\circ\rr_o,\ \ \lim_{n\to\infty} \int_{M}\big\{|f_n- f|^2+|\nn(f_n-f)|^2\big\}\d\theta=0.\end{equation}
Noting that   $\eta\mapsto \eta (\psi\circ\rr_o)$ obeys the Gamma-distribution with parameter
$\dd:= \theta(\psi\circ\rr_o)<\infty$,  we have
$$\int_{\MM} F_0^2\,\d\BG_\theta=\int_{\MM} \eta(\psi\circ\rr_o)^2\, \BG_\theta(\d\eta)= \dd(\dd+1)<\infty.$$ By the dominated convergence theorem and \eqref{*09}, the functions $G_n(\eta):= \eta(f_n), n\ge 1$ satisfy
$$\lim_{n\to\infty} \int_{\MM} |F_0-G_n|^2\d\BG_\theta=0.$$
Moreover,   \eqref{DD}, \eqref{SQF}, \eqref{EENU} and \eqref{*09}  imply
\beg{align*}& \lim_{n,m\to\infty}  \EE_{\BG_\theta}^\ll(G_n-G_m, G_n-G_m)=\lim_{n,m\to\infty} \int_{\MM}\GG^\ll(G_n-F_0,G_n-F_0)\d\BG_\theta\\
&=\lim_{n,m\to\infty} \int_{\MM} \eta(|f_n-f)|^2+\ll |f-f_n|^2)   \BG_\theta(\d\eta)\\
&= \lim_{n,m\to\infty} \int_{M} (|f_n-f)|^2+\ll |f-f_n|^2) \d\theta=0.\end{align*}
Therefore, $F_0\in \D(\EE_{\BG_\theta}^\ll).$ Now, let $F_n$ be defined in \eqref{LLFN}.
Then $F_n\to F$ in $L^2(\BG_\theta)$, and \eqref{DD}, \eqref{SQF} and \eqref{EENU} imply
\beg{align*} &\EE_{\BG_\theta}^\ll(F_n-F,F_n-F) \\
&\le 2 \int_{\MM}\Big\{\|F\|_\infty \GG^\ll((F_0-n)^+\land 1,(F_0-n)^+\land 1)
 + \|\GG^\ll(F,F)\|_\infty   | (F_0-n)^+\land 1|^2\Big\}\d\BG_\theta\\
&\le C\int_{\MM} \Big\{\GG^\ll(F_0,F_0) 1_{\{n\le F_0\le n+1\}} +  1_{\{F_0\ge n\}}\Big\}\d\BG_\theta\to 0\ \text{as}\ n\to\infty,\end{align*}  where
$C:= 2(\|F\|_\infty+ \|\GG^\ll(F,F)\|_\infty)<\infty.$ Moreover, the dominated convergence theorem implies that $\|F_n-F\|_{L^2(\BG_\theta}\to 0$ as $n\to\infty$. Therefore, \eqref{NNN} holds as desired.

(c) Spectral gap estimates. By Theorem \ref{T2.1} and \eqref{*PY3}, for any $F\in \F C_b^\infty(M)$ with $\BG_\theta(F)=0$, we have
$$\EE_{\BG_\theta}^\ll(F,F) = \EE_{\hat\theta}^\GG(F\circ\Phi, F\circ\Phi) \ge \ll \pi_{\hat\theta}(|F\circ\Phi|^2) =\ll \BG_\theta(F^2).$$
This implies $\gap(\EE_{\BG_\theta}^\ll)\ge \ll.$
  On the other hand,
  we take $F_0(\eta)= \eta(M)$ and intend to show that $F_0\in \D(\EE_{\BG_\theta}^\ll)$ with
  \beq\label{SPC2} \ff{\EE_{\BG_\theta}^\ll(F_0,F_0)}{\BG_\theta(F_0^2)- \BG_\theta(F_0)^2}= \ll,\end{equation} so that by definition
  $\gap( \EE_{\BG_\theta}^\ll)\le  \ll,$ and hence the proof is finished.

 Let $\xi\in C_0^\infty([0,\infty))$ such that $0\le \xi\le 1, \xi(s)=1$ for $s\le 1.$ For $h_n$ in \eqref{HNN}, define
\beq\label{SPC3} F_n(\eta)= \int_0^{\eta(h_n)} \xi(s/n)\d s,\ \ \eta\in \MM.\end{equation}    Then $F_n\in \F C_b^\infty(M)$ and $0\le F_n\uparrow F_0$ as $n\uparrow \infty$.
 By \eqref{BDD0} we have
\beq\label{FST}\beg{split} & \int_\MM F_0\d\BG_\theta = \ff{1}{\GG(\theta(M))} \int_0^\infty s^{\theta(M)} \e^{-s}\d s = \theta(M),\\
&\int_\MM F^2_0\d\BG_\theta = \ff{1}{\GG(\theta(M))} \int_0^\infty s^{\theta(M)+1} \e^{-s}\d s = \theta(M)^2+\theta(M).\end{split}\end{equation}
So, the dominated convergence theorem  implies  $\BG_\theta(|F_n-F_0|^2) \to 0$ as $n\to\infty$, and
\beg{align*}  &\lim_{n,m\to\infty} \EE_{\BG_\theta}^\ll (F_n-F_m, F_n-F_m)\\
 &= \lim_{n,m\to\infty} \int_\MM \Big\{\eta\big(\ll\big|\xi(\eta(h_n)/n) h_n- \xi(\eta(h_m)/m) h_m\big|^2\big)\\
  &\qquad\qquad\qquad \qquad  + \big(\big|\xi(\eta(h_n)/n) \nn h_n- \xi(\eta(h_m)/m) \nn h_m\big|^2\big)\Big\}\BG_\theta(\d\eta) =0.\end{align*}
Therefor, $F_0\in \D(\EE_{\BG_\theta}^\ll)$ with
\beg{align*} \EE_{\BG_\theta}^\ll(F_0,F_0)&= \lim_{n\to\infty} \EE_{\BG_\theta}^\ll(F_n,F_n) = \ll \int_\MM \eta\big(\ll\big|\xi(\eta(h_n)/n) h_n\big|^2 + \big|\xi(\eta(h_n)/n) \nn h_n\big|^2\big)\BG_\theta(\d\eta)\\
 &= \ll \BG_\theta (F_0)=\ll\theta(M).\end{align*}   Combining this with \eqref{FST} we derive \eqref{SPC2}
and hence finish the proof.

\end{proof}

\section{Proof of Theorem \ref{T1.1}(2)}
The quasi-regularity can be proved by the same means as in the step (b) of the proof of Theorem \ref{T1.1}(1). So, we need only to prove the closability, the formula \eqref{*DR}, and the claimed spectral gap bounds.

(1) Closability and \eqref{*DR}. Let $F_0(\eta)=\eta(M)$. It suffices to prove that for any $F\in \F C_b^\infty(M)$, one has $F_0\cdot (F\circ\Psi) \in \D(\EE_{\BG_\theta}^\ll)$ such that \eqref{*DR} holds. Indeed, since
 $(\EE^\ll_{\BG_\theta}, \D(\EE^\ll_{\BG_\theta}))$ is closed, \eqref{*DR} implies the closability of $(\EE_{\BD_\theta}^\ll,\F C_b^\infty(M))$.

Similarly to (b) in the proof of Theorem \ref{T3.1}, and noting that $\nn^{int}F_0=0, \nn^{ext} F_0=1$, we see that for any $G\in \F C_b^\infty(M)$,   $GF_n\to GF_0$ in $\D(\EE_{\BG_\theta}^\ll)$ with
\beq\label{ABC1} \EE_{\BG_\theta}^\ll(F_0G,F_0G)= \int_{\MM}\eta\big(\ll \big|F_0\nn^{ext}G+ G\big|^2 +\big|F_0\nn^{int}G\big|^2\big)\BG_\theta(\d\eta).\end{equation}
However, since for  general $F\in \F C_b^\infty(M)$ we do not have $F\circ\Psi\in \F C_b^\infty(M)$, this does not imply $F_0\cdot (F\circ\Psi)\in \D(\EE_{\BG_\theta}^\ll)$ as desired.

To   approximate $F\circ\Psi$ using functions in $\F C_b^\infty(M)$,
we write $F= f(\<h_1,\cdot\>,\cdots,\<h_k,\cdot\>) $ for some $f\in C_b^\infty(\R^k)$ and $h_1,\cdots, h_k\in C_0^\infty(M)$.
Since the Riemannian manifold $M$ is complete, we may construct  $\{\phi_n\}_{n\ge 1}^\infty \subset C_0^\infty(M)$ such that
\beq\label{PPH} 0\le \phi_n\uparrow 1,\ \ |\nn \phi_n|\le \e^{-n}, \ \ \cup_{i=1}^k {\rm supp} h_i\subset \{\phi_n=1\},\ \ n\ge 1.\end{equation}
  For any $n\ge 1,$ let
$$\tt F_n (\eta)= F_0 F_n,\ \ F_n(\eta):= f\Big(\ff{\<h_1,\eta\>}{\eta(\phi_n)+n^{-1}},\cdots, \ff{\<h_k,\eta\>}{\eta(\phi_n)+n^{-1}}\Big),\ \ \eta\in\MM. $$
So, \eqref{ABC1} implies $\tt F_n\in \D(\EE_{\BG_\theta}^\ll)$.
Obviously,  $\BG_\theta(|\tt F_n- F_0(F\circ \Psi)|^2)\to 0\ (n\to\infty)$.  By   \eqref{PPH}  and noting that
\beg{align*}&\nn^{ext} F_n(\eta)= \sum_{i=1}^k (\pp_i f)\Big(\ff{\<h_1,\eta\>}{\eta(\phi_n)+n^{-1}},\cdots, \ff{\<h_k,\eta\>}{\eta(\phi_n)+n^{-1}}\Big)\Big(\ff{h_i} {\eta(\phi_n)+n^{-1}}- \ff{\eta(h_i)\phi_n}{(\eta(\phi_n)+n^{-1})^2}\Big),\\
& \nn^{int} F_n(\eta)= \sum_{i=1}^k (\pp_i f)\Big(\ff{\<h_1,\eta\>}{\eta(\phi_n)+n^{-1}},\cdots, \ff{\<h_k,\eta\>}{\eta(\phi_n)+n^{-1}}\Big)\Big(\ff{\nn h_i} {\eta(\phi_n)+n^{-1}}- \ff{\eta(h_i)\nn \phi_n}{(\eta(\phi_n)+n^{-1})^2}\Big),
\end{align*} we may find out a constant $c>0$ such that
\beg{align*} I_{n,m}(\eta)&:= \eta\Big(\ll \big|\eta(M) \big(\nn^{ext}(F_n-F_m)(\eta)+ (F_n-F_m)(\eta)\big|^2 + \big|\eta(M) \nn^{int}(F_n-F_m)(\eta)\big|^2\Big)\\
&\le c(1+\eta(M)),\ \ \eta\in \MM,\ \ n,m\ge 1.\end{align*} Since $F_0\in L^1(\BG_\theta)$, by
   \eqref{ABC1} and Fatou's lemma, we arrive at
\beg{align*}\limsup_{n,m\to\infty} \EE_{\BG_\theta}^\ll(\tt F_n-\tt F_m, \tt F_n-\tt F_m)&= \limsup_{n,m\to\infty} \int_\MM I_{n,m}(\eta)\BG_\theta(\d\eta)\\
&\le  \int_{\MM} \Big(\limsup_{n,m\to\infty} I_{n,m}(\eta)\Big)\BG_\theta(\d\eta)=0.
\end{align*} Then $F_0(F\circ\Psi)=\lim_{n\to\infty} \tt F_n\in \D(\EE_{\BG_\theta}^\ll)$.  Similarly, for $G=g(\<h_1,\cdot\>,\cdots, \<h_k,\cdot\>)\in \F C_b^\infty(M)$, we define $\tt G_n= F_0 G_n$ in the same way. By the above formulas of intrinsic and extrinsic derivatives for $F_n$ and $G_n$,  it is easy to see that
 $$\lim_{n\to\infty} \nn^{ext} F_n(\eta)= \ff 1 {\eta(M)} (\tt\nn^{ext}F)(\Psi(\eta)),\ \ \lim_{n\to\infty} \nn^{int} F_n(\eta)= \ff 1 {\eta(M)} (\nn^{int}F)(\Psi(\eta))$$ and the same holds for $(G_n,G)$ replacing $(F_n,F)$. Therefore, by the dominated convergence theorem and using
   \eqref{BDD0} and \eqref{BDD},  we obtain
\beg{align*} &\EE_{\BG_ \theta}^\ll (F_0(F\circ\Psi),F_0(G\circ\Psi))=\lim_{n\to\infty} \EE_{\BG_\theta}^\ll(\tt F_n, \tt G_n)\\
&= \lim_{n\to\infty} \int_\MM \Big[\eta\big(\ll\big\{\eta(M) (\nn^{ext}F_n)(\eta)+ F_n(\eta)\big\}\cdot \big\{\eta(M) (\nn^{ext}G_n)(\eta)+ G_n(\eta)\big\}\big) \\
&\qquad\qquad \qquad\qquad \qquad\qquad +\eta(M)^2\eta\big(\big\<(\nn^{int}F_n)(\eta), (\nn^{int}G_n)(\eta)\big\>_{M}\big)\Big]\BG_\theta(\d\eta)\\
&= \int_\MM \eta(M) \Big[ \ll \Psi(\eta)\big(\{ \tt\nn^{ext} F (\Psi(\eta))\}\cdot \{ \tt\nn^{ext} G (\Psi(\eta))\}\big)+ \ll  (FG)(\Psi(\eta))  \\
 &\qquad \qquad\qquad  \qquad\qquad\qquad+ \Psi(\eta)\big(\<\nn^{int} F (\Psi(\eta)),\nn^{int} G (\Psi(\eta))\>_M \big)\Big]\BG_\theta(\d\eta) \\
&= \theta(M) \int_{\MM_1} \eta\Big(\ll \{\tt\nn^{ext} F(\eta)\}\cdot \{\tt\nn^{ext} G(\eta)\} + \ll (FG)(\eta) +\< \nn^{int}F(\eta), \nn^{int}G(\eta)\>_M \Big) \BD_\theta(\d\eta)\\
&= \ll\theta(M)\BD_\theta(FG) +\theta(M) \EE_{\BD_\theta}^\ll(F,G). \end{align*}
Therefore, \eqref{*DR} holds.

\

 (2) Spectral gap estimates. Since
 $$\EE_{\BD_\theta}^\ll(F,F)\ge \ll  \EE_{\BD_\theta}^{FV}(F,F),\ \ F \in \F C_b^\infty(M),$$
 where $\EE_{\BD_\theta}^{FV}$ is given in \eqref{FV} as the Dirichlet form of the Fleming-Viot process, it follows from \eqref{FVG} that
 $$\gap(\EE_{\BD_\theta}^\ll)\ge \ll  \gap(\EE_{\BD_\theta}^{FV})=\ll\theta(M).$$

To verify the spectral gap upper bound, for any $h\in C_b^\infty(M)$ with $\theta(h)=0$ and $\theta(h^2)=1$, let
 $F_h(\eta)= \eta(h).$ Then $F_h\in \D(\EE_{\BD_\theta}^\ll).$ It suffices to show that
\beq\label{KKL}  0<    \EE_{\BD_\theta}^\ll(F_h,F_h)\le \big\{\ll\theta(M)+ \theta(|\nn h|^2)(\theta(M)+1)\}\cdot\{\BD_\theta(F_h^2)- \BD_\theta(F_h)^2\}.\end{equation}
To this end, we recall that (see for instance the proof of \cite[Lemma 7.2]{RW01})
\beq\label{PPP} \pi_{\hat\theta}(\<\hat h, \cdot\>)= \hat\theta(\hat h),\ \ \pi_{\hat\theta}(\<\hat h, \cdot\>^2)= \hat\theta(\hat h^2)+\hat\theta(\hat h)^2,\ \ \hat h\in L^2(\hat\theta).\end{equation} Letting
$\hat g(x,s)=s g(x)$ for $g\in L^2(M,\theta)$ and applying \eqref{BDD0}, \eqref{BDD}, and \eqref{GGT}, we deduce from \eqref{PPP} that
\beg{align*} &\theta(M) \BD_\theta(F_g)= \int_{\MM}\eta(M) \<g,\Psi(\eta)\>\BG_\theta(\d\eta)= \int_\MM \eta(g)\BG_\theta(\d\eta)\\
&= \pi_{\hat\theta}(\<\hat g,\cdot\>)= \theta(g),\ \ g\in L^2(M,\theta),\end{align*}
and similarly,
\beg{align*} &\theta(M)(\theta(M)+1) \BD_\theta(|F_g|^2)= \int_{\MM}\eta(M)^2  \<g,\Psi(\eta)\>^2\BG_\theta(\d\eta)= \int_\MM |\eta(g)|^2\BG_\theta(\d\eta)\\
&= \pi_{\hat\theta}(\<\hat g,\cdot\>^2)= \theta(g^2)+\theta(g)^2,\ \ g\in L^2(M,\theta).\end{align*} Thus, for $h\in C_b^\infty(M)$ with $\theta(h)=0$ and $\theta(h^2)=1$, we have
\beg{align*} & \BD_\theta(F_h^2)- \BD_\theta(F_h)^2= \ff{\theta(h^2)-\theta(h)^2}{\theta(M)(\theta(M)+1)} - \ff{\theta(h)^2}{\theta(M)^2} = \ff 1 {\theta(M)(\theta(M)+1)},\end{align*} and
\beg{align*} &\EE_{\BD_\theta}^\ll(F_h,F_h) =\BD_\theta\big(\ll  \big\{\<h^2,\cdot\>-\<h,\cdot\>^2\big\}+ \<|\nn h|^2,\cdot\>\big)\\
&=\ll  \Big(\ff{\theta(h^2)}{\theta(M)}-\ff{\theta(h^2)+\theta(h)^2}{\theta(M)(\theta(M)+1)}\Big)  + \ff{\theta(|\nn h|^2)}{\theta(M)}
 = \ff{\ll}{\theta(M)+1} + \ff{\theta(|\nn h|^2)}{\theta(M)}.\end{align*} Therefore, \eqref{KKL} holds.

 \paragraph{Acknowledgement.} We would like to thank  the referees for helpful comments on an earlier version of the paper.

\end{document}